\documentclass[10pt, conference]{IEEEtran}
\ifCLASSINFOpdf
  \usepackage[pdftex]{graphicx}
\else
  \usepackage[dvips]{graphicx}
\fi

\usepackage{flushend}
\usepackage{caption} 
\usepackage{amsmath}
\usepackage{amssymb}
\usepackage{color}

\def\N{{\mathbb N}}

\def\F{{\mathbb F}}
\def\E{{\mathbb E}}
\def\L{{\mathbb L}}

\newtheorem{theorem}{Theorem}
\newtheorem{lem}{Lemma}

\newcommand{\vectii}[2]{\ensuremath{\begin{pmatrix} #1 \\ #2 \end{pmatrix}}}

\begin{document}
%
\title{Analytical Results for a Small Multiple-layer Parking System}

\author{\IEEEauthorblockN{Sjoert Fleurke}
\IEEEauthorblockA{Radiocommunications Agency Netherlands\\
Postbus 450, 9700 AL Groningen\\
The Netherlands\\ Email: sjoert.fleurke@agentschaptelecom.nl}
\and
\IEEEauthorblockN{Aernout C.D. van Enter}
\IEEEauthorblockA{Johann Bernoulli Institute for Mathematics\\
and Computer Science, University of Groningen\\
Nijenborgh 9, 9747 AG Groningen\\
The Netherlands\\ Email: a.c.d.van.enter@rug.nl}
}


%


\maketitle

\begin{abstract}
In this article a multilayer parking system of size $n=3$ is studied. We prove that the asymptotic limit of the particle density in the center approaches a maximum of 1/2 in higher layers. This means a significant increase of capacity compared to the first layer where this value is 1/3. This is remarkable because the process is solely driven by randomness. We conjecture that this result applies to all finite parking systems with $n \ge 2$.
\end{abstract}


%
\IEEEpeerreviewmaketitle

\section{Introduction}

\begin{figure}[h]
\setlength{\unitlength}{0.5mm}
\begin{center}
\begin{picture}(80,80)(5,-10)
\multiput(15,12)(19,0){4}{\line(1,0){15}}
\multiput(15,25.1)(19,0){4}{\line(1,0){15}}
\multiput(15,38.2)(19,0){4}{\line(1,0){15}}
\multiput(15,51.3)(19,0){4}{\line(1,0){15}}
\multiput(15,64.4)(19,0){4}{\line(1,0){15}}
\multiput(32,14)(19,0){3}{\line(0,1){9.2}}
\multiput(32,27.2)(19,0){3}{\line(0,1){9.2}}
\multiput(32,40.4)(19,0){3}{\line(0,1){9.2}}
\multiput(32,53.6)(19,0){3}{\line(0,1){9.2}}
\multiput(32,66.9)(19,0){3}{\line(0,1){5}}
\multiput(32,11.7)(19,0){3}{\circle{4.1}}
\multiput(32,25)(19,0){3}{\circle{4.1}}
\multiput(32,38.3)(19,0){3}{\circle{4.1}}
\multiput(32,51.6)(19,0){3}{\circle{4.1}}
\multiput(32,64.9)(19,0){3}{\circle{4.1}}

\multiput(10,10)(77,0){2}{\small{$\times$}}
\put(51,11.7){\circle*{4.1}} \put(32,25.0){\circle*{4.1}}
\put(51,38.3){\circle*{4.1}} \put(51,11.7){\circle*{4.1}}

\put(1,10){\small{1}} \put(1,23){\small{2}} \put(1,36){\small{3}} \put(1,49){\small{4}} \put(1,62){\small{...}}
\put(8,0){\small{-2}} \put(28,0){\small{-1}} \put(49,0){\small{0}} \put(68,0){\small{1}} \put(87,0){\small{2}}

\put(34,54){\small{A}} \put(53,54){\small{B}} \put(72,27){\small{C}}
\put(-5,55){\small{r}}  \put(78,-5){\small{x}}
\put(75,-8){\vector(1,0){11}} \put(-9,52){\vector(0,1){11}}
\end{picture}
\end{center}
\caption{Parking lattice consisting of 3 positions where parking is allowed. Three particles have already arrived consecutively at positions 0, -1, and 0. The next particle will be deposited either in A, B, or C depending on the position (-1, 0, or 1 respectively) where it arrives. The `$\times$' symbols at -2 and 2 indicate that at those x-positions no particles arrive.}
\label{fig:lattice}
\end{figure}
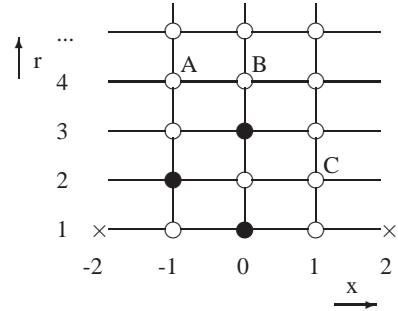

Suppose we have a lattice $\L(x,r)$ consisting of sites $(x,r)$ with positions $x \in \{-2,-1,0,1,2\}$ and heights $r \in \N$. At each position particles arrive according to independent Poisson processes $N_t(x)$. We impose boundary conditions $N_t(-2) = N_t(2) = 0$. The particles pile up across the layers but they are not allowed to ``interfere'' with particles earlier deposited in neighboring sites at the same layer. In other words, the horizontal distance between two particles has to be at least 2. Furthermore, the model has no screening i.e. the particles are always deposited in the lowest possible layer (see Fig. \ref{fig:lattice}).

Our model can be formulated more precisely in the following way.
\begin{enumerate}
  \item The state-space is $\F := (\L, \N^+)^{\{0,1\}}$.
  \item The process $\kappa_t(x,r) = 1$ if there is a particle at $(x,r)$ at time $t$ and 0 otherwise.
  \item When a particle arrives at site $x$ at time $t$, it will be deposited at $h_t(x) := \min\{r : \kappa_t(y,r) = 0, \forall_{y \in N_x} \}$, where neighborhood set $N_x$ consists of site $x$ and the sites with distance 1 from it.
\end{enumerate}
The density $\rho_t(x,r)$ of a site at $(x,r) \in \L$ is defined as the expectation of the occupancy of that site at time $t$, or $\rho_t(x,r) := \E \kappa_t(x,r)$. The end-density of a site is $\rho_{\infty}(x,r)$.

Our models can be viewed either as particle deposition, car parking \cite{Renyi} \cite{Hemmer}, or as models for random sequential adsorption \cite{Evans}.
In this article we will use the terminology of particle deposition.
We focus on the densities of the sites in the center, i.e. those with coordinates $(0,r), r \in \N^+$. The majority of the existing literature in which discrete parking is analytically treated is about monolayer models \cite{Hemmer} \cite{Co62} \cite{Fl08}, while most literature about multi-layer models is based on simulations \cite{Privman} \cite{Dehling07}. However, in \cite{Kulske} it was shown that in an infinite parking system the second layer has a higher capacity than the first layer and in \cite{Fleurke2011} time-dependent density formulas for the first few layers of small finite parking systems are calculated.

In this paper we continue the work on calculating the particle densities in a small multi-layer parking model. We hope our result will lead to further insights also in systems with bigger sizes.

\section{Particle Densities in the Case of Deposition at 3 Vertices on an Interval}
In this section we will analytically calculate the end-densities in the case of a system with 3 vertices.

\begin{theorem}
\label{theo:1}
Consider a multilayer parking system with three vertices. The average density at vertex 0 at height $h+1 \ge 1$ and at time $t$ obeys the following formula
\begin{equation}
\begin{split}
&\rho_t^{(h+1)}(0) = \sum_{k=0} ^h \left[ \vectii{h}{k} \vectii{h+k}{k} \left(\frac{1}{3}\right)^{h+k+1} \right. \cr
&~ \left. + 2\vectii{h}{k} \sum_{j=0} ^{k-1} \vectii{h+j}{j} \left(\frac{1}{3}\right)^{h+j+1} \right] \cr
&~ - \sum_{k=0} ^h \left[ \vectii{h}{k} \vectii{h+k}{k} \left(\frac{1}{3}\right)^{h+k+1} \sum_{i=0} ^{h+k} \frac{(3t)^i}{i!}  \right. \cr
&~ +\!\! \left. 2 \vectii{h}{k} \!\! \sum_{j=0} ^{k-1} \!\! \vectii{h+j}{j} \!\! \left(\frac{1}{3}\right)^{h+j+1} \sum_{i=0}^{h+j} \frac{(3t)^i}{i!} \right]  e^{-3t} \cr
\end{split}
\end{equation}
\end{theorem}
\subsection{Proof of Theorem 1}
The proof of this result is based on the fact that a new particle that arrives at $x$ at time $t$ will always be deposited in
layer $h_t(x)+1$. Therefore the derivative of the density at a height $y+1$ at time $t$ is equal to the probability that $H_t(x) = y$.

For the height stochastic variable $H_t(0)$ we can state that
\begin{lem}
\begin{eqnarray}
H_t(0) &=& N_t(0) + \max(N_t(-1), N_t(1)) 
\end{eqnarray}
where $N_t(x)$ is the number of Poisson arrivals at site $x$ at time
$t$.\\
\end{lem}

\begin{IEEEproof}\\
Recall from the Introduction that the height $H_t(0)$ at position 0 is defined as the total number of layers containing one or two particles. So, we may write
\begin{equation}
\begin{split}
H_t(0) &= \sum_{r=1} ^{\infty} \kappa_t(-1,r) + \kappa_t(0,r) + \kappa_t(1,r) \cr
&~- \kappa_t(-1,r) \kappa_t(1,r) \cr
&= N_t(0) + N_t(-1) + N_t(1) \cr
&~- \sum_{r=1} ^{\infty} \kappa_t(-1,r) \kappa_t(1,r)
\end{split}
\label{form:h0}
\end{equation}
The value of the last term may be written as
\begin{equation}
\begin{split}
&\sum_{r=1} ^{\infty} \kappa_t(-1,r) \kappa_t(1,r) \cr
&~= \left\{
	\begin{array}{ll}
		N_t(-1)  & \mbox{if } N_t(-1) \leq N_t(1) \\
		N_t(1)   & \mbox{if } N_t(-1) > N_t(1)
	\end{array}
\right.
\end{split}
\end{equation}
or more simply
\begin{eqnarray}
\sum_{i=1} ^{\infty} \kappa_t(-1,r) \kappa_t(1,r)  = \min(N_t(-1),N_t(1))
\end{eqnarray}
 Combining this result with (\ref{form:h0}) completes the proof of the lemma.
\end{IEEEproof}

The next step is to calculate the probability $\Pr(H_t(0) = h)$. Therefore we first need to derive the density of the term $\max(N_t(-1), N_t(1))$.
\begin{lem}
\begin{equation}
\begin{split}
\Pr(\max(N_t(-1), N_t(1) = n) &= \left( e^{-t} \frac{t^{n}}{n!}\right)^2 \cr
&~+ 2 e^{-t}\frac{t^{n}}{n!} e^{-t}  \sum_{j=0} ^{n-1} \frac{t^j}{j!}
\end{split}
\end{equation}
\end{lem}
\begin{IEEEproof}
\begin{equation}
\begin{split}
&\Pr \left( \max(N_t(-1), N_t(1)) = n \right) \cr
                   &= \Pr( [N_t(-1)=n] \cap [N_t(1) < n] ) \cr
                   &~+ \Pr( [N_t(1)=n] \cap [N_t(-1) < n]) \cr
                   &~~+ \Pr(N_t(-1)= N_t(1) = n) \cr
                   &= 2\Pr([N_t(-1)=n] \cap [N_t(1) < n] ) \cr
                   &~+ \Pr(N_t(-1)= N_t(1) = n) \cr
                   &= 2\Pr(N_t(-1)= n) \Pr(N_t(1) < n) \cr
                   &~+ \Pr(N_t(1) = n)^2  \cr
                   &= 2 e^{-t}\frac{t^{n}}{n!} e^{-t} \sum_{j=0} ^{n-1}  \frac{t^j}{j!} + \left( e^{-t} \frac{t^{n}}{n!}\right)^2 \cr
                   &= 2 e^{-t}\frac{t^{n}}{n!} e^{-t} \sum_{j=0} ^{n}  \frac{t^j}{j!} - e^{-t} \frac{t^{n}}{n!}
\end{split}
\end{equation}
\end{IEEEproof}

The combination of lemma 1 and lemma 2 provides us a useful expression for the height. Since the probability that $H_t(x) = y$ equals the derivative of the density of the site at height $y+1$ at time $t$ we can continue as follows.

\begin{IEEEproof}
\begin{equation}
\begin{split}
\dot{\rho}_t^{(h+1)}(0) &= \Pr(H_t(0) = h) \cr
                   &= \Pr(N_t(0) + \max(N_t(-1), N_t(1) = h) \cr
                   &= \sum_{k=0} ^h \Pr(N_t(0) = h-k)\cr
                   &~\times \Pr(\max(N_t(-1), N_t(1) = k) \cr
                   &= \sum_{k=0} ^h e^{-t} \frac{t^{h-k}}{(h-k)!} \left( 2 e^{-t}\frac{t^{k}}{k!}  \sum_{j=0} ^{k-1} e^{-t} \frac{t^j}{j!} \right. \cr
                   &~+ \left. \left( e^{-t} \frac{t^{k}}{k!}\right)^2 \right) \cr
&= 2 t^h e^{-3t} \sum_{k=0} ^h \left[ \frac{1}{k!(h-k)!} \sum_{j=0} ^{k-1}  \frac{t^{j}}{j!} \right] \cr
&~ + t^h e^{-3t}  \sum_{k=0} ^h \frac{1}{k!(h-k)!}\frac{t^{k}}{k!}
\end{split}
\end{equation}

Integrating this expression results in the time-dependent densities
$\rho_t^{h+1}$ for layer $h+1$. So, we have

\begin{equation}
\begin{split}
\rho_t^{(h+1)}(0) &= 2 \sum_{k=0} ^h \left[ \frac{1}{k!(h-k)!}
\sum_{j=0} ^{k-1} \frac{\int_0 ^t x^{h+j} e^{-3x} dx }{j!} \right] \cr
&~ + \sum_{k=0} ^h \frac{\int_0 ^t x^{h+k} e^{-3x}  dx  }{(h-k)!k!^2 }%
\end{split}
\end{equation}
Now we use the identity
\begin{equation}
\int e^{-a x} x^{s} dx = -\frac{s!}{a^{s+1}} e^{-ax}\sum_{i=0} ^s \frac{(ax)^i}{i!} \label{form:gamma}
\end{equation}
and get
\begin{equation}
\begin{split}
\rho_t^{(h+1)}(0) &= 2 \sum_{k=0} ^h \left[ \frac{1}{k!(h-k)!} \sum_{j=0} ^{k-1} \frac{(h+j)!}{j!}  \right. \cr
\end{split}
\nonumber
\end{equation}
\begin{equation}
\begin{split}
&~ \times \left. \frac{\left(1 - e^{-3t}\sum_{i=0}^{h+j} \frac{(3t)^i}{i!} \right) }{3^{h+j+1}}\right] \cr%
&~~ + \sum_{k=0} ^h \frac{(h+k)!}{k!^2(h-k)!} \frac{\left( 1 - e^{-3t}\sum_{i=0} ^{h+k} \frac{(3t)^i}{i!} \right)}{3^{h+k+1}} \cr
&= 2 \sum_{k=0} ^h \left[ \vectii{h}{k} \sum_{j=0} ^{k-1} \vectii{h+j}{j} \right. \cr
&~ \times \left. \frac{\left(1 - e^{-3t}\sum_{i=0}^{h+j} \frac{(3t)^i}{i!} \right) }{3^{h+j+1}}\right] \cr%
&~~ + \sum_{k=0} ^h \vectii{h}{k} \vectii{h+k}{k} \frac{\left( 1 - e^{-3t}\sum_{i=0} ^{h+k} \frac{(3t)^i}{i!} \right)}{3^{h+k+1}} %
\end{split}
\end{equation}
\begin{equation}
\begin{split}
&= \sum_{k=0} ^h \left[ \vectii{h}{k} \vectii{h+k}{k} \frac{1}{3^{h+k+1}} \right.~~~~~~~~~~~~~~~~\! \cr
&~ + 2\vectii{h}{k} \left. \sum_{j=0} ^{k-1} \vectii{h+j}{j} \frac{1}{3^{h+j+1}} \right] \cr
&~ - \sum_{k=0} ^h \left[ \vectii{h}{k} \vectii{h+k}{k} \frac{\sum_{i=0} ^{h+k} \frac{(3t)^i}{i!}   }{3^{h+k+1}} \right. \cr
&~   \left. + 2 \vectii{h}{k} \sum_{j=0} ^{k-1} \vectii{h+j}{j} \frac{\sum_{i=0}^{h+j} \frac{(3t)^i}{i!}}{3^{h+j+1}} \right]  e^{-3t}
\nonumber
\end{split}
\end{equation}
This may be rewritten as (with $r = h+1$).
\begin{equation}
\begin{split}
\rho_t^{(r)}(0) &= \sum_{k=0} ^{r-1} \left[ \frac{\vectii{r-1}{k} \vectii{r+k-1}{k}}{3^{r+k}} \right. \cr
&~ \left.+ 2\vectii{r-1}{k} \sum_{j=0} ^{k-1} \vectii{r+j-1}{j} \frac{1}{3^{r+j}} \right] \cr
&~ - \sum_{k=0} ^{r-1} \left[ \vectii{r-1}{k} \vectii{r+k-1}{k} \frac{\sum_{i=0} ^{r+k-1} \frac{(3t)^i}{i!}   }{3^{r+k}} \right. \cr
&~ \left. + 2 \vectii{r-1}{k} \sum_{j=0} ^{k-1} \vectii{r+j-1}{j} \right. \cr
&~~ \times \left. \frac{\sum_{i=0}^{r+j-1} \frac{(3t)^i}{i!}}{3^{r+j}} \right]  e^{-3t}\cr
\end{split}
\end{equation}
\end{IEEEproof}
\begin{figure}[!b]
\begin{picture}(200,140)
     \put(-35,0){\includegraphics[scale=0.40]{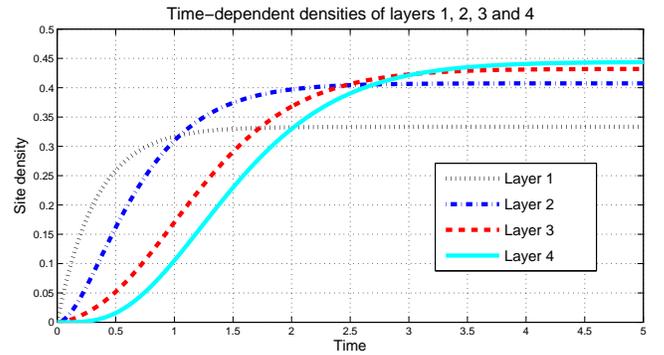}}
\end{picture}
\caption{Particle densities at the sites $(0,r)$ as a function of time in the cases $r$ is 1, 2, 3, and 4 according to (\ref{form:results}).}
\label{fig:timedependent}
\end{figure}
For the first few layers Theorem \ref{theo:1} provides:
\begin{equation}\begin{split}
\rho^{(1)}_t &= \frac{1}{3} - \frac{1}{3} e^{-3 t} \cr
\rho^{(2)}_t &= \frac{11}{27} - \left( \frac{11}{27} + \frac{11}{9}  t + \frac{1}{3} t^2 \right)e^{-3 t} \cr
\rho^{(3)}_t &= \frac{35}{81} - \left( \frac{35}{81} + \frac{35}{27} t + \frac{35}{18} t^2 + \frac{7}{9} t^3 + \frac{1}{12} t^4 \right)e^{-3 t} \cr
\rho^{(4)}_t &= \frac{971}{2187} - \left( \frac{971}{2187} + \frac{971}{729} t + \frac{971}{486} t^2 + \frac{971}{486} t^3\right. \cr
&~ \left. + \frac{283}{324} t^4 + \frac{17}{108} t^5 \right)e^{-3 t} \cr
\end{split}
\label{form:results}
\end{equation}
Confer \cite{Fleurke2011} where the first 3 layers were calculated using a different approach.
%
%
%
%
\section{Calculation of the End-densities}

Close inspection of (\ref{form:results}) reveals that as time goes to infinity the densities of the first 4 layers tend towards $\frac{1}{3}$, $\frac{11}{27}$,  $\frac{35}{81}$, and $\frac{971}{2187}$ respectively. Calculating end-densities for higher layers can be done directly from Theorem \ref{theo:1}. See Table \ref{tab:enddens} for the exact values of the end-densities of the first 10 layers and its decimal approximations.

\begin{table}[h]
\centering
\begin{tabular}{|l|c|c|}
\hline
Layer & End-density & Approximately \\
\hline
1 & $\frac{1}{3}$                   & 0.3333    \\ [1ex]
2 & $\frac{11}{27}$                 & 0.4074    \\ [1ex]
3 & $\frac{35}{81}$                 & 0.4321    \\ [1ex]
4 & $\frac{971}{2187}$              & 0.4440    \\ [1ex]
5 & $\frac{8881}{19683}$            & 0.4512    \\ [1ex]
6 & $\frac{80811}{177147}$          & 0.4562    \\ [1ex]
7 & $\frac{733209}{1594323}$        & 0.4599    \\ [1ex]
8 & $\frac{6640491}{14348907}$      & 0.4628    \\ [1ex]
9 & $\frac{60067809}{129140163}$    & 0.4651    \\ [1ex]
10 & $\frac{542880971}{1162261467}$ & 0.4671    \\ [1ex]
\hline
\end{tabular}
\caption{End-densities calculated using Theorem \ref{theo:1}}
\label{tab:enddens}
\end{table}

\begin{figure}[!h]
\begin{picture}(200,140)
     \put(-30,0){\includegraphics[scale=0.39]{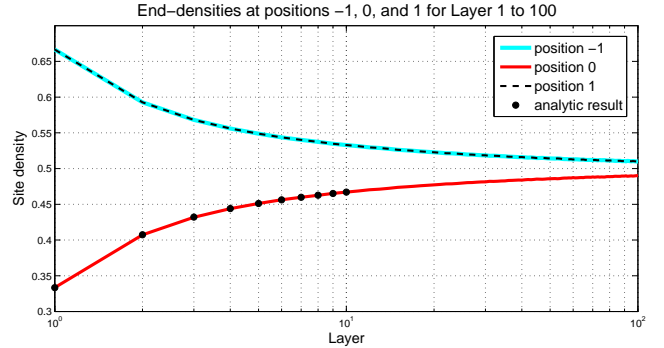}}
\end{picture}
\caption{End-densities as a function of the layer created by $10^7$ simulations. The analytic results from Table \ref{tab:enddens} are plotted as well to demonstrate its concurrence.}
\label{fig:enddensities}
\end{figure}

A plot of these constants for the first 100 layers is shown in Figure \ref{fig:enddensities}. It can be seen clearly that the graph of these end-densities appears to approach the value of $\frac{1}{2}$. In this section we will prove that this is indeed the case.

Define $\rho^{(r)} := \lim_{t \rightarrow \infty} \rho_t(0,r)$. Then we have the following result.

\begin{theorem}
The density at high layers converges in time to the value
\begin{equation}
\lim_{r \rightarrow \infty}  \rho^{(r)} = \frac{1}{2}
\end{equation}
\end{theorem}
To prove this we can take the result of Theorem \ref{theo:1} and focus on the constant term. %
\begin{equation}
\begin{split}
\rho^{(h+1)} &= \sum_{k=0} ^h \left[ \vectii{h}{k} \vectii{h+k}{k} \left(\frac{1}{3}\right)^{h+k+1} \right. \cr
&~ \left. + 2\vectii{h}{k} \sum_{j=0} ^{k-1} \vectii{h+j}{j} \left(\frac{1}{3}\right)^{h+j+1} \right]
\end{split}
\end{equation}
We may rewrite this more conveniently as
\begin{equation}
\begin{split}
\label{firstandsecondterm}
\rho^{(r)} &= \frac{1}{2} \sum_{k=0} ^{r-1} \Pr(X_{r-1,\frac{1}{2}} = k) \Pr(Y_{r,\frac{1}{3}} = k) \cr
&+ \sum_{k=0} ^{r-1}  \Pr(X_{r-1,\frac{1}{2}} = k) \sum_{j=0} ^{r-k} \Pr(Y_{r,\frac{1}{3}} = j)
\end{split}
\end{equation}
Where we used the notation $r = h+1 \in \N^+$, $X_{n,p} \sim B(n,p)$ or $\Pr(X_{n,p} = k) = \vectii{n}{k} (1-p)^{n-k} p^k$, and also $Y_{r,p} \sim NB(r,p)$ or $\Pr(Y_{r,p} = k) = \vectii{r+k-1}{k} (1-p)^r p^k$. We will treat the first and second term of (\ref{firstandsecondterm}) separately in the following lemmas.
\begin{lem}
The first term of (\ref{firstandsecondterm}) converges to zero when $r \rightarrow \infty$, or
\begin{equation}
\begin{split}
\lim_{r \rightarrow \infty} \frac{1}{2} \sum_{k=0} ^{r-1} \Pr(X_{r-1,\frac{1}{2}} = k) \Pr(Y_{r,\frac{1}{3}} = k) = 0
\end{split}
\end{equation}
\end{lem}

\begin{IEEEproof}
\begin{equation}
\begin{split}
&\sum_{k=0} ^{r-1} \Pr(X_{r-1,\frac{1}{2}} = k) \Pr(Y_{r,\frac{1}{3}} = k) \cr
&~=\sum_{k=0} ^{r-1} \Pr(X_{r-1,\frac{1}{2}} = k \cap Y_{r,\frac{1}{3}} = k) \cr
&~= \sum_{k=0} ^{r-1} \Pr(Y_{r,\frac{1}{3}} = X_{r-1,\frac{1}{2}}| X_{r,\frac{1}{3}} = k) \Pr(X_{r,\frac{1}{3}} = k) \cr
&~= \Pr(Y_{r,\frac{1}{3}} = X_{r-1,\frac{1}{2}})
\end{split}
\end{equation}
This represents the probability that the number of successes (with $\Pr(Success) = \frac{1}{3}$) after $r$ failures equals the number of successes in a Binomial experiment of $r-1$ trials and $\Pr(Success) = \frac{1}{2}$. When we let $r \rightarrow \infty$ both $X_{r-1,\frac{1}{2}}$ and $Y_{r,\frac{1}{3}}$ will converge to continuous Gaussian distributions, so that this probability vanishes.
\end{IEEEproof}

\begin{lem} The second term of (\ref{firstandsecondterm}) converges to $1/2$, or
\begin{equation}
\begin{split}
\lim_{r \rightarrow \infty} \sum_{k=0} ^{r-1}  \Pr(X_{r-1,\frac{1}{2}} = k) \sum_{j=0} ^{r-k} \Pr(Y_{r,\frac{1}{3}} = j) = \frac{1}{2} \cr%
\end{split}
\end{equation}
\end{lem}
\begin{IEEEproof}
\begin{equation}
\begin{split}
&\sum_{k=0} ^{r-1}  \Pr(X_{r-1,\frac{1}{2}} = k) \sum_{j=0} ^{r-k} \Pr(Y_{r,\frac{1}{3}} = j) \cr
    &= \sum_{k=0} ^{r-1}  \Pr(X_{r-1,\frac{1}{2}} = k) \Pr(Y_{r,\frac{1}{3}} \le r-k)
\end{split}
\end{equation}
Now we will use the symmetry of the negative binomial distribution for large $r$. Note that $\Pr(Y_{r,p} < r-k) = \Pr(Y_{r,p} > k)$ in this case where $p=1/3$.
\begin{equation}
\begin{split}
&\lim_{r \rightarrow \infty} \frac{1}{2} \left( \sum_{k=0} ^{r-1} \Pr(X_{r-1,\frac{1}{2}} = k) \Pr(Y_{r,\frac{1}{3}} < r - k) \right. \cr
&~  \left. + \sum_{k=0} ^{r-1} \Pr(X_{r-1,\frac{1}{2}} = k) \Pr(Y_{r,\frac{1}{3}} < k+1  \right)\cr
&~= \lim_{r \rightarrow \infty} \frac{1}{2} \left( \sum_{k=0} ^{r-1} \Pr(X_{r-1,\frac{1}{2}} = k) \Pr(Y_{r,\frac{1}{3}} > k) \right. \cr
&~+ \left.  \sum_{k=0} ^{r-1} \Pr(X_{r-1,\frac{1}{2}} = k) \Pr(Y_{r,\frac{1}{3}} < k+1)  \right)\cr
                                &\approx \lim_{r \rightarrow \infty} \frac{1}{2}\sum_{k=0} ^{r-1} \Pr(X_{r-1,\frac{1}{2}} = k) = \frac{1}{2}
\end{split}
\end{equation}
\end{IEEEproof}

\subsection{Alternative proof}

We note however that this result also follows from the following consideration. After a while the differences in height between position -1 and 1 increase to the order of the square root of the total number of dropped particles. This follows by application of the Central Limit Theorem to $K_t := |N_t(-1) - N_t(1)|$.\\
The probability that a new particle drops at a side vertex happens with probability 2/3. So, the probability that this particle raises the height equals 1/2 times 2/3, which is 1/3. This equals the probability that a particle drops on the center vertex, by which the height always increases. For 1/3 of the dropped particles the height does not increase. Thus half of the newly filled layers contains an occupied center vertex, and half will contain two occupied side vertices, which implies density one half.

\subsection{Larger parking systems}
\begin{figure}[!h]
\begin{picture}(100,140)
     \put(-35,0){\includegraphics[scale=0.40]{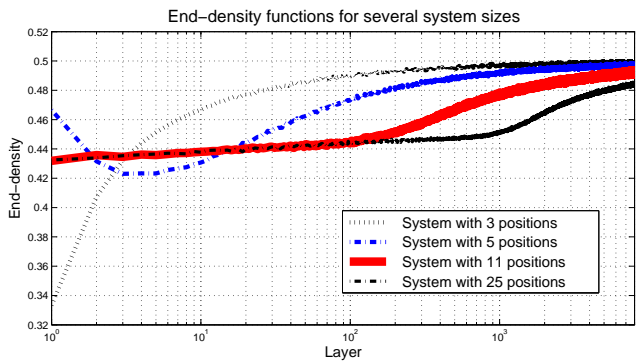}}
\end{picture}
\caption{The fact that the end-densities of the sites in the center converge to 1/2 in the case of 3 positions is not unique. Simulation results of bigger systems suggest that this behavior is not uncommon for finite-sized systems. However, it appears that the bigger the system, the more layers it takes to approach the limit of 1/2.}
\label{fig:enddensitiessimulations}
\end{figure}

The calculation of (end-)densities in larger systems is much more complicated. It is always possible
to calculate the densities on the first layer \cite{Co62} or the first few layers \cite{Fleurke2011} but
going beyond the first few layers in systems with bigger sizes probably requires more advanced methods.\\
However, it is interesting to ask oneself whether the behavior of the small system demonstrated in this article
does also appear in larger systems. Do systems of bigger sizes also generally have higher end-densities in higher layers than
in lower layers, and if so, do those end-densities ultimately approach the maximum value of $\frac{1}{2}$ as well?\\
We conjecture that this is the case for all finite-sized systems. Although we are not able to give hard evidence for
this we can provide some simulation results (Figure \ref{fig:enddensitiessimulations}) supporting our view and justifying further research.

\section{Conclusion}
In this paper we introduced a parking system consisting of 3 positions. The formula for the time dependent densities of the center position for all layers was analytically derived.
Although similar work has been done on the model with screening (see \cite{Sudbury12}, \cite{FlKu10}) to our knowledge this is the first time that densities in a multi-layer particle deposition model without screening were calculated analytically.\\
We paid special attention to the densities of the center sites when $t \rightarrow \infty$, the so called end-densities. We proved that they increase as a function of the layer number and eventually approach the density $1/2$.\\
We showed that in the case of a small system with 3 positions it can be easily understood why the end-density converges to this value. But this is not the case with larger systems although our simulation results do suggest similar end-density behavior. Although not yet fully understood, it thus seems that these randomly driven finite parking systems tend to use the parking space of the center positions more efficiently over time.


\section*{Acknowledgment}
This research was partially supported by the Radiocommunications Agency of the Netherlands.


\begin{thebibliography}{1}

\bibitem{Renyi}
A.~R\'{e}nyi, On a One-dimensional Problem Concerning Random
Space-filling, {\em Publ. Math. Inst. Hung. Acad. Sci.}~{\bf 3},
(1958), pp.~109--127.
\bibitem{Hemmer}
P.C.~Hemmer, The Random Parking Problem, {\em J. Stat. Phys.}~{\bf
57}, (1989), pp.~865--869.

\bibitem{Evans} J.W.~Evans, Random and Cooperative Sequential Adsorption,
{\em Rev. Mod. Phys.}~{\bf 64}(4), (1993), pp.~1281--1327.



\bibitem{Co62} R.~Cohen, H.~Reiss, Kinetics of Reactant Isolation I. One-Dimensional
Problems, {\em J. Chem. Phys.}~{\bf 38}(3), (1963), pp.~680--691.

\bibitem{Fl08} H.G.~Dehling, S.R.~Fleurke, C.~K\"{u}lske,  Parking on a Random Tree,
{\em J. Stat. Phys.}~{\bf 133}(1), (2008), pp.~151--157.
\bibitem{Privman} 
P.~Nielaba, V.~Privman, Multilayer Adsorption with Increasing Layer
Coverage. {\em Phys. Rev.}~{\bf A 45}, (1992), pp.~6099--6102.
\bibitem{Dehling07}
H.G.~Dehling, S.R.~Fleurke, The Sequential
Frequency Assignment Process, {\em Proc. of the 12th WSEAS Internat.
Conf. on Appl. Math.} Cairo, Egypt, (2007), pp.~280--285
%
\bibitem{Kulske}
S.R.~Fleurke, C.~K\"{u}lske, A Second-row Parking Paradox,
{\em J. Stat. Phys.}~\textbf{136}(2), (2009), pp.~285--295.
%
\bibitem{Fleurke2011}
S.R.~Fleurke, {\em Multilayer Particle Deposition Models}, (Groningen Thesis), published by VDM Verlag Dr. Muller, Saarbr\"{u}cken, (2011), p.~33.
\bibitem{Sudbury12}
T.S.~Mountford, A.~Sudbury, Deposition processes with Hardcore Behavior,
{\em J. Stat. Phys.}~{\bf 146}, (2012), pp.~687--700.
\bibitem{FlKu10}
S.R.~Fleurke, C.~K\"ulske,
Multilayer Parking with Screening on a Random Tree, {\em J. Stat.
Phys.}~\textbf{139}(3), (2010), pp.~417--431.

\end{thebibliography}
\end{document}